\documentclass[a4paper]{jpconf}
\usepackage{graphicx}
\usepackage[latin1]{inputenc}
\usepackage{amsmath}
\usepackage{amsfonts}
\usepackage{rotating}
\usepackage{amssymb}

\newtheorem{teorema}{Teorema}[section]

\newtheorem{remark}[teorema]{Remark}
\newtheorem{prop}[teorema]{Proposition}

\newcommand{\n}{\nabla}
\renewcommand{\a}{\alpha}
\renewcommand{\b}{\beta}
\newcommand{\oln}{\overline{\nabla}}
\newcommand{\ol}{\overline}
\newcommand{\mc}{\mathcal}

\def\M{{\ol{\mc{M}}}}

\begin{document}
\title{Second-order symmetric Lorentzian  manifolds II: structure and global properties}

\author{O F Blanco$^{1}$, M S\'anchez$^{1}$ and J M M Senovilla$^2$}

\address{$^{1}$ Departamento de Geometr\'ia y Topolog\'ia , Facultad de Ciencias, Universidad de Granada
Campus Fuentenueva s/n, 18071 Granada, Spain
}
\address{$^2$ F\'isica Te\'orica, Universidad del Pa\'is Vasco, Apartado 644, 48080 Bilbao, Spain}

\ead{oihane@ugr.es, sanchezm@ugr.es, josemm.senovilla@ehu.es}

\begin{abstract}
We give a summary of recent results on the explicit local form of the second-order symmetric 
Lorentzian manifolds in arbitrary dimension, and its global version. These spacetimes turn out to be essentially a
specific subclass of plane waves.

\end{abstract}

\section{Introduction}

Last year we presented the local explicit form of the {\it second-order symmetric} 
(or {\it $2$nd symmetric}\footnote{The notion of {\em
3-symmetric space} (without ordinal), which is widely spread in the literature, was introduced by Gray \cite{GR} in a different
context.} for short) spacetimes in four dimensions \cite{O1}. The techniques used for that purpose were specific of the dimension. 
Our aim now is to introduce new techniques that generalize the result to arbitrary dimension. 
Recall that $2$nd symmetric spacetimes are characterized by the condition 
$$\n_\rho \n_\sigma R^\alpha~_{\beta\lambda\mu}=0.$$
They were introduced by one of the authors in \cite{SN1} as the simplest natural generalization of the {\it locally symmetric} 
spacetimes (those satisfying $\n_\sigma R^\alpha~_{\beta\lambda\mu}=0$). 
In the Riemannian \cite{Be},\cite{CA},\cite{He} and Lorentzian \cite{CW} cases, the locally symmetric spaces are well known. 
However, $2$nd symmetric spaces have been hardly considered in the literature, probably because the local de Rham decomposition 
implies that $2$nd symmetry and local symmetry are equivalent in the proper Riemannian case \cite{NO},\cite{TA}. 
In the Lorentzian case such decomposition fails \cite{WU3}, and so does the equivalence. 
Thus, even though in arbitrary signature {\it semi-symmetric} spaces 
(those satisfying $\n_{[\rho} \n_{\sigma ]} R^\alpha~_{\beta\lambda\mu}=0$, see for example \cite{SZ1},\cite{SZ2} 
for the Riemannian case and \cite{IE}, \cite{SV} for the Lorentzian one) have always been regarded as the simple natural 
generalization of the locally symmetric ones, in indefinite signatures the $2$nd symmetry is a simpler natural choice. 

In this contribution we announce the following classification of {\it proper $2$nd symmetric} spacetimes 
(here proper means that they are not locally symmetric, that is $\n_\sigma R^\alpha~_{\beta\lambda\mu}\neq 0$) 
in arbitrary dimension \cite{O2}:

\begin{teorema}\label{lem4}
 A  $n$-dimensional proper $2$nd symmetric spacetime $(M,g)$ is
 locally isometric 
 to a product spacetime $(M_1\times M_2,g_1\oplus
 g_2)$ where
 $(M_2,g_2)$ is a non-flat Riemannian symmetric space and $(M_1,g_1)$ is a generalized Cahen-Wallach
 space  of order $2$,
 defined as $M_1=\mathbb{R}^{d+2}$ ($d\geq 1$) endowed with the metric
$$
g_1=-2du\left(dv+p_{ij}(u)x^ix^j
du\right)+\delta_{ij}dx^i dx^j ,~~i,j=2,\ldots,d+1
$$
\noindent where $\{u,v,x^i\}$ are the Cartesian
coordinates of $\mathbb{R}^{d+2}$, each function $p_{ij}$ is affine, i.e., $p_{ij}(u)=\a_{ij} u+\b_{ij}$ for some $\a_{ij}, \b_{ij}\in \mathbb{R}$, and $\sum_{i,j=2}^{d+1}(\a_{ij})^2\neq 0$.

If in addition $(M, g)$ is geodesically complete and simply
connected, then $(M, g)$ is globally isometric to one of such
products.


\end{teorema}

\section{Sketch of the proof of theorem \ref{lem4}}

As in the four-dimensional case, the starting point for the proof is \cite[Theorem 4.2]{SN1}, which states when, 
in the Lorentzian case, the $2$nd symmetry condition does not imply local symmetry. This can only happen
when the spacetime possesses a null (non-vanishing) covariantly constant vector field. Therefore, it is locally isometric to a 
Brinkmann spacetime, whose metric can be written in an appropiate {\it Brinkmann chart} $\{u,v,x^i\}$ as \cite{BR},\cite{ZA}:
$$
ds^2=-2 du (dv+H du+W_i dx^i)+g_{ij}dx^idx^j,\hspace{0.2cm} i,j\in\{2,\ldots,n-1\},
$$\noindent where $H$, $W_i$ and $g_{ij}=g_{ji}$ are functions independent of $v$, otherwise arbitrary. With this data at hand, the steps of the proof are the following:
\begin{enumerate}
 \item \label{1} A suitable choice of vector field basis simplifies the computation of the equations of $2$nd symmetry.
\item \label{2} As a first consequence, the leaves of the foliation $\M$ generated by the submanifolds with constant 
values of $u$ and $v$ are locally symmetric Riemannian manifolds for any Brinkmann chart.
\item Using the integrability equations and (\ref{2}), one proves that the only non-vanishing components of the curvature in 
the basis (\ref{1}) are $\n_0 R^1~_{i0j}$.
\item Applying a generalization of the Eisenhart theorem \cite{Eis}, ``the spatial part'' $g_{ij}$ is splits
into two parts, one of them is Ricci-flat which, on using a result in \cite{Ale}, becomes actually flat.
 \item  After lengthy computations, one can prove that there exists a Brinkmann chart such that the metric becomes:
\begin{eqnarray}\nonumber ds^2&=&-2 du \left(dv+\left(H_1(u,x^a)+H_2(u,x^{a'})\right) du+
W_a(u,x^b) dx^a+W_{a'}(u,x^{b'}) dx^{a'}\right)\\\nonumber &+&\delta_{ab}dx^adx^b+g_{a'b'}dx^{a'}dx^{b'}, 
\end{eqnarray} 
with  $a,b\in\{2,\ldots,d+1\},\, a',b'\in\{d+2,\ldots,n-1\},\, g_{a'b'}\neq \delta_{a'b'}.$
\item At this stage, the equations can be reorganized in two blocks: 
\begin{itemize}
 \item one corresponding to the equations of local symmetry, whose solutions are already known \cite{CW} (see also \cite{Be}),
 \item the other one corresponding to the equations of $2$nd symmetry for a plane wave, which are easily solvable.
\end{itemize}
\end{enumerate}

To develop the proof, the introduction of some mathematical tools is required, namely:
\begin{enumerate}
 \item[(1)] Some operators defined on the tensor bundle $T^r_s \M$, denoted as $D_0$ (transverse operator) and $\oln$ (which arises from the covariant derivative on the leaves), and an adaptation of the exterior derivative to $\M$.
\item[(2)] The curvature tensor $\ol{R}^s~_{ijk}$ associated intrinsically to the foliation (i.e., to $\oln$), as well as its
naturally associated Ricci tensor $\ol{R}_{ij}$ and scalar curvature $\ol{R}$.
\item[(3)] Some tensor fields on $\M$ defined as projections of $R^\alpha~_{\beta\lambda\mu}$ and 
$\n_{\rho}R^\alpha~_{\beta\lambda\mu}$, as well as some algebraic lemmas on vector spaces with positive definite inner product.
\item[(4)] An appropriate adaptation of the Eisenhart theorem to Brinkmann spacetimes, as mentioned above. 
This involves the operators $D_0$ and $\overline{\nabla}$ and ensures the reducibility of the entire family of metrics $g_{ij}$ 
---which depend on $u$--- simultaneously.
\end{enumerate}

Before dealing with the global part of the proof, we  define 
\begin{itemize}
 \item {\it proper $r$-th symmetric} spacetimes as those satisfying $\n_{\rho_1}\ldots\n_{\rho_r} R^\alpha~_{\beta\lambda\mu}=0$ but 
$\n_{\rho_1}\ldots\n_{\rho_{r-1}} R^\alpha~_{\beta\lambda\mu}\neq 0$,
\item the generalized {\it Cahen-Wallach spacetimes of order $r$}, denoted by $CW^d_r(A)$, as the spacetimes 
$(\mathbb{R}^d,g_A)$, $d\geq 3$ with metric
$$
g_A=-2du\left(dv+ A_{ij}x^ix^j du\right)+
\delta_{ij}dx^i dx^j,
$$ 
where  $A=(A_{ij})$ is a square symmetric matrix of dimension $d-2$ such that 
$A\equiv A(u)=A^{(r-1)}u^{r-1}+\ldots+A^{(1)} u+ A^{(0)}, A^{(r-1)}\neq 0$, 
and each $A^{(l)}$ is a square symmetric matrix of constants, $l=1,\ldots,r-1$. 
Observe that these spacetimes are proper $r$th symmetric \cite{SN1},
as well as analytic and geodesically complete \cite{CFS}. 
\end{itemize}
Then, the global result follows directly from:
\begin{prop}
 For any matrix $A$ as above:
\begin{itemize}
 \item the direct product of $CW^n_{r}(A)$ and a finite number of Riemannian symmetric 
(non-necessarily simply connected) spaces is a proper $r$-th symmetric Lorentzian manifold, 
which is in fact geodesically complete and analytic, and
\item if a complete Lorentzian manifold $(M,g)$ is locally isometric to spaces of the form given above, then its universal covering is also of this type.
\end{itemize}
 \end{prop}

\begin{remark}{\rm 
The techniques used to solve the arbitrary dimensional case differ from the techniques used in four dimensions 
since in the latter case the Petrov classification of the Weyl tensor \cite{ES} was used. 
Such classification is not available in arbitrary dimension with the same simplicity, 
though some generalizations have already apeared in the literature  (see \cite{Co1},\cite{Co3},\cite{SN2}). 
Therefore, we had to refocuse the problem in arbitrary dimension using a different, and more powerful, viewpoint.}
 
\end{remark}


\ack
The three authors are partially supported by the grant P09-FQM-4496 (J. Andaluc\'ia), OFB and MS are partially 
supported by MTM2007-60731 (Spanish MICINN with FEDER funds), and JMMS is supported by grants FIS2010-15492 (MICINN) 
and GIU06/37 (UPV/EHU).

\end{document}